\documentclass[12pt]{amsart}
\usepackage{amsfonts}
\usepackage{amssymb}
\usepackage{latexsym}

\newtheorem{theorem}{Theorem}[section]
\newtheorem{lemma}[theorem]{Lemma}
\newtheorem{proposition}[theorem]{Proposition}
\newtheorem{corollary}[theorem]{Corollary}

\numberwithin{equation}{section}

\def\qed{$\hfill\square$}
\def\im{\mbox{Im}\,}
\def\proof{{\it Proof. }}

\begin{document}
\title{Injective Envelopes of Separable C$^*$-algebras}

\thanks{This research is supported in part by the Natural Sciences and Engineering Research Council of Canada}
\thanks{ 2000 Mathematics Subject Classification:  Primary 46L05;
Secondary 46L07 }
\thanks{Keywords and Phrases:
\emph{injective envelope, local multiplier algebra, regular monotone completion, C$^*$-algebra,
AW$^*$-algebra}}
\thanks{June 15, 2005}

\author[Mart\'in Argerami]{{Mart\'\i n~Argerami}}
\author[Douglas R.~Farenick]{{Douglas R.~Farenick}}
\address
        {Department of Mathematics and Statistics\\
    University of Regina \\
    Regina, Saskatchewan\\
    Canada
        \ S4S 0A2}
\email{argerami@math.uregina.ca}
\email{farenick@math.uregina.ca}

\begin{abstract}
Characterisations of those separable C$^*$-algebras that have type
I injective envelopes or W$^*$-algebra injective envelopes are
presented.
\end{abstract}

\maketitle
\markboth{\textsc{M.~Argerami and D.R.~Farenick}}
{\textsc{Injective Envelopes of Separable C$^*$-algebras}}

An operator system $I$ is injective if for every inclusion $E\subset F$ of operator systems
each completely positive linear map $\omega:E\rightarrow I$ has a completely positive
extension to $F$.
An injective envelope of an operator system $E$ is an injective operator system $I$ such that $E\subseteq I$
and $I$ is minimal among all injective operator systems that contain $E$. That is,
if $E\subseteq I_0\subseteq I$, with $I_0$ injective, then $I_0=I$.
Hamana \cite{hamana1979} proved that every operator system $E$ has an injective envelope and that all injective
envelopes of $E$ are completely isometric.
Because every injective operator system is completely order isomorphic
to an injective C$^*$-algebra \cite{choi--effros1977}, and because two C$^*$-algebras are
$*$-isomorphic if and only if they are
completely order isomorphic \cite{choi1974}, one can unambiguously refer to ``the" injective envelope of $E$,
which is an injective C$^*$-algebra $I(E)$ that contains $E$ as an operator system.

If $E$ is a C$^*$-algebra $A$, then $A$ is contained in $I(A)$ as a C$^*$-subalgebra.
The purpose of the present paper is to study how
properties of a C$^*$-algebra $A$ determine properties of its injective envelope,
especially in the case of separable C$^*$-algebras $A$.

The injective envelope $I(A)$ of any C$^*$-algebra $A$ is a monotone complete C$^*$-algebra. Thus, $I(A)$
is a direct sum of AW$^*$-algebras of types I, II, and III. Herein we show that if $A$ is separable,
then $I(A)$ has no direct summand that is finite and of type II.
Further, we show that a separable C$^*$-algebra $A$ has
a type I injective envelope if and only if $A$ has a liminal essential ideal. We also characterise
those separable C$^*$-algebras $A$ for which $I(A)$ is a W$^*$-algebra.

There are a number of other useful enveloping structures that contain a given C$^*$-algebra $A$ as a C$^*$-subalgebra.
Of these, the local multiplier algebra $M_{\rm loc}(A)$
\cite{ara--mathieu2003,frank--paulsen,somerset2000} and the
regular monotone completion $\overline A$ \cite{wright1976a,ozawa--saito,hamana1981}
have important roles in arriving at our
our results. These structures, together with the injective envelope, are discussed in the following
preliminary section.

\section{Preliminary Results}

\subsection{Terminology and notation}
As usual, we will denote by $B(H)$ and $K(H)$ the set of bounded
and compact operators on a Hilbert space $H$. Because the algebras
under study are not represented in any particular way as acting on
a Hilbert space, we shall employ the following terminology. A
C$^*$-algebra $B$ is said to be a $W^*$-algebra if, as a Banach
space, $B$ is the dual space $X^*$ of some (in fact, unique)
Banach space $X$. It is a classical fact \cite[Theorem III.3.5]{tak}
that a C$^*$-algebra $B$ is a W$^*$-algebra if and only if $B$ has
a representation as a von Neumann algebra of operators acting on
some Hilbert space. A C$^*$-algebra $B$ is called an
AW$^*$-algebra if the left annihilator of each right ideal in $A$
is of the form $Ap$ for some projection $p\in A$. Although every
W$^*$-algebra is an AW$^*$-algebra, the converse is not true:
there exist AW$^*$-algebras that fail that have any faithful
representation as a von Neumann algebra.

If $B$ is an AW$^*$-algebra, then $p\sim q$ denotes the Murray--von Neumann
equivalence of projections $p$ and $q$ in $B$. Thus, a projection $p\in B$ is finite
if $q\sim p$ and $q\le p$ only if $q=p$; otherwise $p$ is an infinite projection. If
the identity $1\in B$ is a finite projection, then
$B$ is said to be finite algebra.
A projection $p\in B$ is abelian if the AW$^*$-algebra $pBp$ is commutative.

An AW$^*$-algebra $B$ is said to be: of type I if every
direct summand of $B$ has an abelian projection; of type II if $B$ has no abelian
projections but every direct summand has a finite projection; and of type III if
all projections in $B$ are infinite. If the centre $\mathcal Z(B)$ of an AW$^*$-algebra
is $\mathbb C$, then $B$ is a factor.
Type I AW$^*$-algebras are of considerable interest herein. In particular, type I AW$^*$-algebras
are injective C$^*$-algebras \cite{hamana1981} and type I AW$^*$-factors are of the form
$B(H)$ \cite{kaplansky1952}.

AW$^*$-algebras differ from W$^*$-algebras in that the former can fail to have any normal states.
An AW$^*$-algebra $B$ is wild \cite{wright1976b} if the only normal positive linear functional $\varphi$ on $B$
is $\varphi=0$. Every AW$^*$-factor is either a W$^*$-algebra or a wild AW$^*$-algebra \cite{wright1976b}.

A C$^*$-algebra $A$ is said to be postliminal (or type I, or GCR)
if every representation of $A$ generates a type I von Neumann algebra, and $A$ is liminal (or CCR) if every irreducible
representation $\pi:A\rightarrow B(H)$ satisfies $\pi(A)=K(H)$.
An elementary C$^*$-algebra is one that $*$-isomorphic to $K(H)$ for some Hilbert space $H$.

We shall employ the following notation from \cite{ara--mathieu2003}.
If $\{E_\alpha\}_{\alpha\in\Lambda}$ is a family of operator systems, then
\[
\begin{array}{rcl}
\displaystyle\prod_{\alpha\in\Lambda}E_\alpha &=&\left\{(e_\alpha)_{\alpha}\,:\,e_\alpha\in E_\alpha
\mbox{ and }\sup_\alpha\|e_\alpha\|<\infty\right\}\,; \\
\displaystyle\bigoplus_{\alpha\in\Lambda}E_\alpha &=&\{(e_\alpha)_\alpha\,:\,e_\alpha\in E_\alpha
\mbox{ and }
\forall\,\varepsilon>0\\
&&\mbox{ only finitely many }e_\alpha\mbox{ satisfy }\|e_\alpha\|>\varepsilon\}\,.
\end{array}
\]
Note that if $\{A_\alpha\}_{\alpha\in\Lambda}$ is a family of
C$^*$-algebras, then $\prod_{\alpha}A_\alpha$ and $\bigoplus_{\alpha}A_\alpha$ are C$^*$-algebras
and
$\bigoplus_{\alpha}A_\alpha$ is an ideal of $\prod_{\alpha}A_\alpha$. As operator systems can always be realised
as $*$-closed, unital subspaces of unital C$^*$-algebras,
$\prod_{\alpha\in\Lambda}E_\alpha$ is an operator system for every
family $\{E_\alpha\}_{\alpha\in\Lambda}$ of operator systems.

\subsection{Injective envelopes}
An operator system $I$ is injective if for every inclusion $E\subset F$ of operator systems
each completely positive linear map $\omega:E\rightarrow I$ has a completely positive
extension to $F$. Arveson's extension theorem \cite{arveson1969} for completely positive linear maps
with values in $B(H)$
demonstrates that $B(H)$ is injective. This fact can be used to show that
if an operator system $I$ is represented as a unital, $*$-closed subspace
of $B(H)$, then $I$ is injective if and only if $I$ is the range of some completely positive
linear map $\phi:B(H)\rightarrow B(H)$ for which $\phi^2=\phi$. Such maps $\phi$ are commonly referred
to as projections, or conditional expectations. A theorem of Choi and Effros \cite{choi--effros1977}
demonstrates that if $I$ is an injective operator system given by the range of a projection $\phi$ on
$B(H)$, then $I$ is completely order isomorphic to a C$^*$-algebra,
obtained by changing the product of $I$ to $x\circ y=\phi(xy)$.

An injective envelope of an operator system $E$ is an injective operator system $I$ and a
complete isometry $\kappa:E\rightarrow I$ such that, if $I_0$ is an injective operator
system with $\kappa(E)\subseteq I_0\subseteq I$, then $I_0=I$.
The existence and uniqueness (up
to complete isometry) of the
injective envelope was established by
Hamana \cite{hamana1979}; thus, it is a common practice to drop reference to $\kappa$ and assume that $E$ is already
realised as an operator system in $I$.
The following proposition of Hamana is a useful criterion for determining when an injective operator system $I$
containing $E$ is an injective envelope.

\begin{proposition}\label{hamana37} {\rm (\cite[Lemma 3.7]{hamana1979})}
Consider an inclusion $E\subseteq I$ of operator systems, where $I$
is injective. The following statements are equivalent.
\begin{enumerate}
\item $I$ is an injective envelope of $E$.
\item The only completely positive linear map $\omega:I\rightarrow I$ for which $\omega|_E=id_E$
is the identity map $\omega=id_I$.
\end{enumerate}
\end{proposition}

We note below a property that we shall make frequent use of.

\begin{lemma}\label{direct sum of injectives}
If $\{E_\alpha\}_{\alpha\in\Lambda}$ is a family of operator systems, then $\prod_{\alpha}E_\alpha$ is
injective if and only if $E_\alpha$ is injective for every $\alpha\in\Lambda$.
\end{lemma}

\proof
Fix an inclusion $E\subset F$ of operator systems.

Assume that $\prod_\alpha E_\alpha$ is injective. If $\varphi:E\rightarrow E_\beta$ is completely positive,
define $\tilde\varphi:E\rightarrow\prod_\alpha E_\alpha$
by $(\tilde\varphi(x))_\beta=\varphi(x)$ and $(\tilde\varphi(x))_\alpha=0$ if $\alpha\ne\beta$.
Then there exists $\psi:F\rightarrow\prod_\alpha E_\alpha$ extending $\tilde\varphi$. So
$\pi_\beta\circ\psi$ is a completely positive extension of $\varphi$.

Conversely, if $E_\alpha$ is injective for every $\alpha$, and $\varphi:E\rightarrow\prod_\alpha E_\alpha$ is
completely positive, then for each $\alpha$ the map $\pi_\alpha\circ\varphi:E\rightarrow E_\alpha$ is
completely positive, and so there exists $\psi_\alpha:F\rightarrow E_\alpha$ completely positive extension.
Thus the map $\prod_\alpha\psi_\alpha:F\rightarrow\prod_\alpha E_\alpha$
is a completely positive extension of $\varphi$.
\qed\medskip

\subsection{Regular monotone completions}
A C$^*$-algebra $B$ is monotone complete if every
bounded increasing net $\{h_\alpha\}_\alpha$ in $B_{\rm sa}$ has a least
upper bound in $B_{sa}$,
where $B_{\rm sa}$ denotes the real vector space of hermitian elements of $B$.
The least upper bound of a bounded increasing net $\{h_\alpha\}_\alpha$ in $B_{\rm sa}$
is denoted by $\sup_\alpha h_\alpha$.
A C$^*$-algebra $B$ is monotone $\sigma$-complete if every
bounded increasing sequence $\{h_n\}_{n\in\mathbb N}$ in $B_{\rm sa}$ has a least
upper bound in $B_{sa}$. (The terminology ``monotone complete" is called ``monotone
closed" in some of the standard texts, such as \cite{pedersen1979} and \cite{tak}. We follow
Hamana \cite{hamana1981} by using the term ``monotone closed" in a sense different from
\cite{pedersen1979} and \cite{tak}; this
is explained below.)

Monotone complete C$^*$-algebras are unital \cite{tak} and if
$B$ is monotone $\sigma$-complete and satisfies the countable chain condition (namely, for each
for each $S\subset B_{\rm sa}$ that is bounded above in $B_{\rm sa}$ there is a
countable subset $S_0\subseteq S$ such that any upper bound for $S_0$ is also
an upper bound for $S$), then $B$ is monotone complete \cite{wright1976a}.
Every W$^*$-algebra
is monotone complete and
a C$^*$-algebra $B$
is an AW$^*$-algebra if and only if each maximal abelian C$^*$-subalgebra $D\subseteq B$ is monotone
complete. However, it is not known whether every AW$^*$-algebra is monotone complete.
A well-known theorem of Tomiyama \cite{tomiyama} for conditional expectations between C$^*$-algebras,
which is proved below for operator systems, implies in particular that the injective envelope of an operator
system is monotone closed.

\begin{proposition}\label{range}
Let $E\subseteq M$ be operator systems, with $M$ monotone complete. If there exists a positive linear map
$\phi:M\rightarrow E$ such that $\phi_E=id_E$, then $E$ is monotone complete.
\end{proposition}
\proof Let $\{h_\alpha\}_\alpha$ be a bounded increasing net in $E$.
It is in particular an increasing bounded net in $M$,
so there exists $\tilde{h}\in M$, $\tilde{h}=\sup_\alpha h_\alpha$. Let $h=\phi(\tilde{h})$. Then $h-h_\alpha=
\phi(\tilde{h}-h_\alpha)\ge0$, for every $\alpha$, so that $h$ is an upper bound for $\{h_\alpha\}_\alpha$.
If $k\in E$ and $h_\alpha\le k$ for every $\alpha$, then because $k\in M$ we have that
$\tilde{h}\le k$. Thus, $k-h=\phi(k-\tilde{h})\ge0$, which implies that $h$ is the
supremum of $\{h_\alpha\}_\alpha$ in $E$.
\qed\medskip

\begin{corollary} The injective envelope $I(A)$ of any C$^*$-algebra $A$ is monotone complete.
In particular, $I(A)$ is an AW$^*$-algebra.
\end{corollary}

If $B$ is a monotone complete C$^*$-algebra, then a subset $S\subseteq B_{\rm sa}$ is
monotone closed in $B$ if, for every bounded increasing net $\{s_\alpha\}_\alpha$ in $S$,
$\sup_\alpha s_\alpha$ (which exists in $B$) is contained
in $S$. In particular, if $A$ is a C$^*$-subalgebra of $B$ and if
$\mbox{m-cl}_B\,A_{\rm sa}$ denotes the smallest subset of $B_{\rm sa}$ that contains $A_{\rm sa}$
and is monotone closed in $B$, then the monotone closure of $A$ in $B$ is defined to be
the set
\[
\mbox{m-cl}_B\,A\;=\;\mbox{m-cl}_B\,A_{\rm sa}\;+\;i\ \mbox{m-cl}_B\,A_{\rm sa}\,.
\]
It so happens that $\mbox{m-cl}_B\,A$ is a monotone complete C$^*$-subalgebra of
$B$ \cite[Lemma 1.4]{hamana1981}.

A C$^*$-subalgebra $C$ of $B$ is called a monotone closed C$^*$-subalgebra of $B$ if
$\mbox{m-cl}_B\,C=C$. Because the property of $C$ being monotone closed in $B$ involves both $C$ and $B$,
it is possible for a C$^*$-subalgebra $C$ of $B$ to be monotone complete yet fail to be
monotone closed in $B$. In fact, it is frequently the case that a von Neumann algebra $M\subset B(H)$
is not monotone closed in $B(H)$.

A C$^*$-subalgebra $A$ of a C$^*$-algebra $B$ is said to be order dense in $B$ if
\[
h\;=\;\sup\{k\in A^+\,:\,k\le h\}\,,\quad\forall\,h\in B^+\,.
\]
For example, $K(H)$ is order dense in $B(H)$.

A regular monotone completion of a C$^*$-algebra $A$ is a C$^*$-algebra $B$ such that
\begin{enumerate}
\item $A$ is a C$^*$-subalgebra of $B$,
\item $B$ is monotone complete,
\item $\mbox{m-cl}_B\,A=B$, and
\item $A$ is order dense in $B$.
\end{enumerate}

In \cite{hamana1981}, Hamana proved that a regular monotone completion exists for every C$^*$-algebra $A$
and any two regular monotone completions of $A$ are $*$-isomorphic. Thus,
$\overline A$ is used to denote ``the" regular monotone completion of $A$.
Hamana's construction of $\overline A$ is via the injective envelope of $A$. Namely, $\overline A$ is the
monotone closure of $A$ in $I(A)$.

The regular monotone $\sigma$-completion $\overline A\,{}^\sigma$ of a C$^*$-algebra $A$
was introduced by Wright \cite{wright1976a}.
Hamana recovers $\overline A\,{}^\sigma$ via the injective envelope by considering
monotone $\sigma$-closure of $A$ in $I(A)$ (the definitions are analogous to earlier ones, but
with sequences in place of nets).

For each C$^*$-algebra $A$ there is a representation in which
\[
A\,\subseteq\, \overline A\,{}^\sigma \,\subseteq\, \overline A\,\subseteq\, I(A)\,,
\]
where each containment is as a C$^*$-subalgebra. We shall assume this representation in our work herein.
An important feature of this sequence of containments is:
\[
\overline A\,\mbox{ is monotone closed in }\,I(A)\,.
\]

\begin{theorem}\label{sigma completions} Assume that $A$ is a separable C$^*$-algebra.
\begin{enumerate}
\item {\rm (Wright)} \cite{wright1976a} $\overline A\,{}^\sigma=\overline A$.
\item {\rm (Ozawa--Sait\^o)} \cite{ozawa--saito}
The AW$^*$-algebra $\overline A$ has no type II direct summand.
\item {\rm (Hamana)} \cite{hamana1981} If $A$ is postliminal, then $\overline A$ is of type I.
\item\label{inj-four} {\rm (Sait\^o)} \cite{saito1980}
If $K\subseteq A$ is an essential ideal of $A$, then $\overline K=\overline A$.
\item\label{inj-five} If $K\subseteq A$ is an essential ideal of $A$, then then $I(K)=I(A)$.
\end{enumerate}
\end{theorem}

\proof Only the proof of (\ref{inj-five}) need be given, as it is not explicitly stated in the literature.
By (\ref{inj-four}), $\overline K=\overline A$ if $K\subseteq A$ is an essential ideal of $A$. Furthermore,
$I(\overline A)=I(A)$, by \cite[Lemma 3.7]{hamana1981}. Hence, $I(K)=I(\overline K)=I(\overline A)=I(A)$.
\qed\medskip

\subsection{Local multiplier algebras}

The multiplier algebra of a C$^*$-algebra $A$ is the C$^*$-subalgebra $M(A)$ of the enveloping von Neumann
algebra $A^{**}$ that consists of all $x\in A^{**}$ for which $xa\in A$ and $ax\in A$, for all $a\in A$.
If $J\subseteq A$ is an ideal, then $J^{**}$ is identified with the closure of $J$ in $A^{**}$ with respect
to the strong operator topology. Thus, if $J$ and $K$ are ideals of $A$, and if $J\subseteq K$, then
$M(J)\supseteq M(K)\supseteq M(A)$.

An ideal $K$ of $A$ is said to be essential if $K\cap J\not=\{0\}$
for every nonzero ideal $J\subseteq A$. Any essential ideal is necessarily nonzero.
Consider the multiplier algebra $M(J)$ of any essential ideal $J$ of $A$. If
$\mathcal E(A)$ is the set of essential ideals of $A$, partially ordered by reverse inclusion, then the
set $\mathfrak E(A)$ of multiplier algebras $M(K)$ of $K\in\mathcal E(A)$ is a directed system of
C$^*$-algebras.
$M_{\rm loc}(A)$ is then defined to be the C$^*$ direct limit of the directed system $K\in\mathcal E(A)$.
In \cite{ara--mathieu2003},
Ara and Mathieu give a systematic account of the theory of local multiplier algebras of C$^*$-algebras.
Their book is our basic reference on the topic.

There are various ways to realise $M_{\rm loc}(A)$ ``concretely'' as a C$^*$-subalgebra
of some other C$^*$-algebra:
\begin{enumerate}
\item[(i)] as a C$^*$-subalgebra of a quotient of $A^{**}$ \cite{ara--mathieu2003};
\item[(ii)] as a C$^*$-subalgebra of a quotient of $A^{**}$,
where the quotient is monotone $\sigma$-complete \cite{somerset2000};
\item[(iii)] as a C$^*$-subalgebra of $I(A)$ \cite{frank--paulsen}.
\end{enumerate}
In this final case, $M_{\rm loc}(A)$ is realised by idealisers in $I(A)$ of essential ideals of $A$. Specifically,
by \cite[Corollary 4.3]{frank--paulsen},
\[
M_{\rm loc}(A)\;=\;\left( \bigcup_{K\in\mathcal E(A)} \{x\in I(A)\,:\,xK+Kx\subseteq K\}
\right)^{-}\,,
\]
where the closure is with respect to the norm topology of $I(A)$. Thus,
\[
A\,\subseteq\,M_{\rm loc}(A) \,\subseteq\, I(A)
\]
is an inclusion of C$^*$-subalgebras. In \cite{frank2002}, Frank showed an additional sequence of
inclusions as C$^*$-subalgebras:
\[
A\,\subseteq\,M_{\rm loc}(A)
\,\subseteq\, M_{\rm loc}\left(M_{\rm loc}(A)\right) \,\subseteq\,\overline A \,\subseteq\, I(A)\,.
\]

\subsection{Injective envelopes of separable and prime C$^*$-algebras}

\begin{proposition}\label{no type II_1} If $A$ is a separable C$^*$-algebra, then $I(A)$
does not have a finite type II direct summand.
\end{proposition}

\proof It is enough to show that if $I(A)$ has a finite direct summand, then this summand is of type I.
Because $I(A)e=I(Ae)$ for any central projection $e\in I(A)$
\cite[Lemma 6.2]{hamana1981}, and since the C$^*$-algebra
$Ae$ is separable, we may assume without loss of generality that $I(A)$ itself is a finite algebra.
Thus, the identity $1\in I(A)$ is a finite projection, and so $1$ is a finite projection in $\overline A$ as well.
Therefore, $\overline A$ is of type I \cite[Theorem 2]{ozawa--saito}.
But type I algebras are injective; hence $\overline A=I(A)$.
\qed\medskip

The next proposition, which builds on work of Hamana, determines which
C$^*$-algebras lead to factors.

\begin{proposition}\label{factors} The following statements are equivalent for any
C$^*$-algebra $A$.
\begin{enumerate}
\item\label{fac1} $\overline A$ is a factor.
\item\label{fac2} $I(A)$ is a factor.
\item\label{fac3} $A$ is prime.
\end{enumerate}
\end{proposition}

\proof The equivalence of (\ref{fac1}) and (\ref{fac3}) was established by Hamana
\cite[Theorem 7.1]{hamana1981}. To prove that (\ref{fac1}) and (\ref{fac2}) are equivalent, note that
$\mathcal Z(\overline A)=\mathcal Z\left( I(\overline A)\right)$,
because $\overline A$ is monotone complete \cite[Theorem 6.3]{hamana1981}.
Further, because $I(\overline A)=I(A)$ \cite[Lemma 3.7]{hamana1981}, we conclude that
$\mathcal Z(\overline A)=\mathcal Z\left( I(A)\right)$. Thus, $\overline A$ is a factor
if and only if $I(A)$ is a factor.
\qed\medskip

\section{W$^*$-algebra Injective Envelopes}

The injective envelope $I(A)$ of any C$^*$-algebra $A$
is an AW$^*$-algebra. However, in rare
instances $I(A)$ is known in fact to be a W$^*$-algebra. This is so if $A$ can be represented as acting on a
Hilbert space in such a way
as to contain every compact operator \cite{arveson1972,hamana1979}.
In this section we characterise those separable C$^*$-algebras $A$ for which $I(A)$ is a W$^*$-algebra.

\begin{lemma}\label{kadison-1956} If $A$ is a C$^*$-algebra for which $I(A)$ is
a W$^*$-algebra, then $\overline A$ is a W$^*$-algebra.
\end{lemma}

\proof Without loss of generality we may assume that $I(A)$ is represented as a von
Neumann algebra acting on a Hilbert space. Let $\{h_\alpha\}_\alpha$ be any bounded increasing
net in $\overline A_{\rm sa}$. Because $I(A)$ is a von Neumann algebra,
$\{h_\alpha\}_\alpha$ has a least upper bound $h$ such that $h=\lim_\alpha h_\alpha$
in the strong-operator-topology. Note that the supremum of $\{h_\alpha\}_\alpha$ in $I(A)$
necessarily coincides with $h$ and, because $\overline A$ is monotone closed in $I(A)$, $h\in\overline A$.
Thus, $\overline A$ is a C$^*$-algebra of operators for which the strong-operator limit of
every bounded increasing net of hermitian elements of $\overline A$ belongs to $\overline A$.
By \cite[Lemma 1]{kadison1956}, this implies that $\overline A$ is a von Neumann algebra.
\qed\medskip

\begin{lemma}\label{charact-minimal}
The following statements are equivalent for a von Neumann algebra $M$.
\begin{enumerate}
\item\label{prop-typeI} $M$ is a direct product of type I factors.
\item\label{prop-minimal} $M$ is generated by its minimal projections.
\end{enumerate}
\end{lemma}

The lemma above is well known. However, as it is important for our work, the ideas that underlie
the proof are worth noting here briefly.
First of all, if
$M$ is a direct product of type $I$ factors, then $M$ is generated by the family of all the minimal
projections of all the factors.
Conversely, if $M$ is generated by minimal projections,
then it cannot have a type II nor type III direct summand. Indeed, if
$Me$ is type II or type III,
with $e$ a central projection in $M$, consider $q\in M$ a minimal projection such that $qe\ne0$. Such
a projection exists because otherwise $e=0$. Since $q$ is minimal in $M$, $qe=qeq=q$ and so $q\in Me$. But then
$Me$ admits a minimal projection, which is a contradiction.
Thus $M$ is type I, and it can be expressed as a direct integral over a type I factor-valued measure. The diffuse
part of this measure has to be zero, because any projection in the diffuse part will not be minimal, and we can
reason as before. Therefore,
the measure is atomic and $M$ is a direct product of type I factors.

\begin{corollary}\label{atomic-injective}
If $M$ is a von Neumann algebra generated by minimal projections, then $M$ is injective.
\end{corollary}
\proof
Type I factors are injective, by Arveson's theorem \cite{arveson1969}
on the injectivity of $B(H)$.
Lemma \ref{charact-minimal} asserts that $M$ is a direct product of type I factors; by
Lemma \ref{direct sum of injectives},
every direct
product of injective C$^*$-algebras is injective. Hence, $M$ is injective.
\qed\medskip

\begin{lemma}\label{prop-sup}
Suppose that $A$ is a C$^*$-subalgebra of a von Neumann algebra $M$ and that $M=A''$.
\begin{enumerate}
\item\label{A has min projs} If $M$ is generated by its minimal projections,
each of which is contained in $A$, then $A$ is order
dense in $M$.
\item\label{sep sup} If $A$ is separable and if $A$ is order dense in $M$, then
$M$ is generated by its minimal projections, each of which is contained in $A$.
\end{enumerate}
\end{lemma}
\proof
For the proof of (\ref{A has min projs}), choose a nonzero $h\in M^+$ and
consider the set
\[{\mathcal F}=\{\,(k_i)\subset A^+:\,\sum_{\mbox{finite}}k_i\le h\}.\]
There is
a strictly positive $\lambda$ in the spectrum $\sigma(h)$ of $h$.
Let $e\in M$ be the spectral
projection $e=e^{h}\left([\lambda,\infty)\right)$, where $e^{h}$ denotes the spectral resolution of $h$.
Thus, $0\not= \lambda e\le he$. Moreover, $e$ majorises a minimal projection $p$ of $M$; by hypothesis,
$p\in A$. Thus, $0\not=\lambda p=e(\lambda p)e\le e(\lambda)e=\lambda e
\le he\le h$, and so $\lambda p\in\mathcal F$, which proves that
${\mathcal F}\ne\emptyset$.
It is clear that ${\mathcal F}$ is inductive under inclusions of those families and so,
by Zorn's Lemma, ${\mathcal F}$ has a maximal family $W$. Since every finite sum of
this family is less than $h$,
\[y=\sup\left\{\sum_{k\in K}k:K\mbox{ is finite and }K\subset W\right\}\le h.\]
If $y\ne h$, then $h-y\in M^+$, and so by the first paragraph there exists
$k\in A^+$ such that $k\le h-y$. If it were true that $k\in W$, then for each
net $(h_i)$ of those finite sums of
elements in $W$ such that $h_i\nearrow y$, the net $(h_i+k)\nearrow y+k$, which contradicts the fact that $y$
is the supremum. Hence, $k\not\in W$. But if $k\not\in W$, then the family $W$ is
not maximal, which is again a contradiction.
Therefore, it must be that $y=h$, which proves that $A$ is order dense in $M$.

For the proof of (\ref{sep sup}), note that because $A$ is separable and $A''=M$, to prove
that $M$ is generated by its minimal projections, each of which is contained in $A$, it is
enough, by \cite[p.~139]{tak}, to prove that any normal state $\omega\in M_*$ is faithful precisely when
its restriction $\omega|_A$ to $A$ is faithful. Thus, let $\omega$ be a normal state on $M$ that is faithful on $A$.
Assume that $\omega(h)=0$, where $h\in M^+$. Because $h=\sup\{k\in A^+: k\le h\}$, we have that
$0\le\omega(k)\le\omega(h)=0$ for each $k\in A^+$ with $k\le h$. Thus, $\omega(k)=0$, which implies that $k=0$
because $\omega$ is faithful on $A$. Hence, $h=0$ and so $\omega$ is faithful on $M$.
\qed\medskip

The following theorem is the main result of the present section.

\begin{theorem}\label{envelope is type I}
The following statements are equivalent for a separable C$^*$-algebra $A$.
\begin{enumerate}
\item\label{un} $I(A)$ is a W$^*$-algebra.
\item\label{deux} $I(A)$ is a discrete type I W$^*$-algebra.
\item\label{trois} There exists a faithful representation $\pi:A\rightarrow B(H)$
such that the von Neumann algebra $\pi(A)''$ is generated by its minimal projections, each of which is contained in
$\pi(A)$.
\item\label{quatre} There exists an ideal $K$ of $A$ such that
\begin{enumerate}
\item\label{quatre-2} $K$ is a minimal essential ideal and
\item\label{quatre-1} $K\cong_*\bigoplus_nK(H_n)$, for some sequence of Hilbert
spaces $H_n$.
\end{enumerate}
\end{enumerate}
\end{theorem}

\proof Assume that $I(A)$ is a W$^*$-algebra. Then there is a faithful representation
$\tilde\pi:I(A)\rightarrow B(H)$ such that $\tilde\pi(I(A))$ is a von Neumann algebra and
$\pi(A)$ is a C$^*$-subalgebra of $\tilde\pi(I(A))$, where $\pi=\tilde\pi_{|A}$. Without loss
of generality, we assume that $I(A)$ is a von Neumann algebra acting on a Hilbert space.
Consider
the regular monotone completion $\overline A$ of $A$, which can be realised as the monotone closure of
$A$ in $I(A)$ by Hamana's theorem \cite[Theorem 3.1]{hamana1981}. Furthermore, because $I(A)$ is a von Neumann algebra,
$\overline A$ is a von Neumann algebra, by Lemma \ref{kadison-1956}. Thus,
$A''\subseteq\overline A''=\overline A$. As $A$ is separable and order dense in $A''$, the von Neumann algebra
$A''$ is generated by its minimal projections, each of which is contained in $A$ (Lemma \ref{prop-sup}).
Furthermore,
by Lemma \ref{charact-minimal}, $A''$ is a direct product of type I factors, which implies that $A''$ is injective
by Corollary \ref{atomic-injective}.
Because $A\subseteq A''\subseteq I(A)$, we conclude that
$A''=\overline A=I(A)$, by minimality of the injective envelope.
This proves that (\ref{un}) $\Rightarrow$ (\ref{deux}) $\Rightarrow$ (\ref{trois}).

We next show that (\ref{trois}) $\Rightarrow$ (\ref{quatre}).
Assume
there exists a faithful representation $\pi:A\rightarrow B(H)$
such that the von Neumann algebra $\pi(A)''$ is generated by its minimal projections, each of which is contained in
$\pi(A)$. Without loss of generality, assume that $A$ is already represented as a subalgebra of $B(H)$ and that $M=A''$
is generated by its minimal projections, each of which lie in $A$.

Let $K\subseteq A$ be the ideal of $A$ generated by the minimal projections of $M$.
We first show that $K$ is an essential ideal, minimal
among all essential ideals of $A$. Suppose that $J\subseteq A$ is a nonzero ideal.
Choose any nonzero $h\in J^+$. As shown in the proof of (\ref{A has min projs}) of
Lemma \ref{prop-sup}, there is a $\lambda>0$
and a spectral projection $e\in M$ of $h$ such that $\lambda e\le he$, and
there is a minimal projection
$p$ of $M$ such that $ep=pe=p$ and $0\not=\lambda p\le php\in J\cap K$. That is, $J\cap K\not=\{0\}$, which
proves that $K$ is an essential ideal of $A$.

Because $M=A''$ is generated by its minimal projections, $M$ is a discrete type
$I$ von Neumann algebra, by Lemma \ref{charact-minimal}.
Hence, there is a faithful normal $*$-representation $\varrho$ of $M$ on a Hilbert
space $H$ of the form $H=\bigoplus_n H_n$ such that $\varrho(K)\subseteq\varrho(A)\subseteq\varrho(M)=\prod_n B(H_n)$.
Obviously, the minimal projections of any $B(H_n)$ are minimal projections of $\varrho(M)$ and are, hence, elements
of $\varrho(K)$. On the other hand,
if $e$ is a minimal projection of $\prod_n B(H_n)$, then $e\in B(H_n)$ for some $n\in\mathbb N$ (for otherwise
$e$ is cut by some minimal central projection).
Therefore, $\bigoplus_n K(H_n)\subseteq\varrho(K)$. However,
$\varrho(K)$ is the smallest C$^*$-algebra that contains the minimal projections of $\varrho(M)$; hence
$\varrho(K)=\bigoplus_n K(H_n)$.
Since $K\simeq_*\bigoplus_n K(H_n)$, it is a minimal essential ideal.

We now prove that (\ref{quatre}) $\Rightarrow$ (\ref{un}).
Suppose that $A$ has a minimal essential
ideal $K$ such that $K\cong_*\bigoplus_n K(H_n)$. Therefore, by \cite[Lemma 1.2.1]{ara--mathieu2003},
\[
M(K)\;=\;M\left(\bigoplus_n K(H_n)\right)\;=\;\prod_n M\left(K(H_n)\right)\;=\;\prod_n B(H_n)\,,
\]
which shows that $M(K)$ is a (type I) W$^*$-algebra.
Furthermore, because $K$ is a minimal essential ideal of $A$,
$M(K)=M_{\rm loc}(A)$ by \cite[Remark 2.3.7]{ara--mathieu2003}. Hence, $M_{\rm loc}(A)$ is an injective
W$^*$-algebra. However, $A\subseteq M_{\rm loc}(A)\subseteq I(A)$ as C$^*$-subalgebras, and so by definition
of the injective envelope, it must be that $M_{\rm loc}(A)= I(A)$, which proves that $I(A)$ is a W$^*$-algebra.
\qed\medskip

\section{Type I Injective Envelopes}

One extension of Arveson's fundamental theorem \cite{arveson1969} on the injectivity of $B(H)$ is a result of
Hamana \cite[Proposition 5.2]{hamana1981} that states that every type I AW$^*$-algebra is injective. The following
theorem describes those separable C$^*$-algebras that have type I injective envelopes.

\begin{theorem}\label{aw-typeI} If $A$ is a separable C$^*$-algebra $A$, then
$I(A)$ is a type I AW$^*$-algebra if and only if $A$ has a liminal essential ideal. If this is
the case, then $\overline A=I(A)$.
\end{theorem}

\proof
Assume that $A$ is separable and has a liminal essential ideal $K$.
Because $\overline A$ and $\overline K$ are isomorphic \cite[Corollary 2.1]{saito1980}
and because $K$ is liminal, $\overline A$
is a type I AW$^*$-algebra \cite[Theorem 6.6]{hamana1981}.
Hence, $\overline A=I(A)$ and
$I(A)$ is a type I AW$^*$-algebra.

Conversely, assume that $I(A)$ is a type I AW$^*$-algebra. Because $\overline A\subseteq I(A)$ and because
$\overline A$ and $I(A)$ have the same type I direct summands \cite[Corollary 6.5]{hamana1981}, we conclude
that $\overline A=I(A)$. Thus, $A$ is order dense in $I(A)$.

Because $I(A)$ is of type I, the C$^*$-subalgebra $J\subset I(A)$ generated by the abelian projections of $I(A)$
is a liminal ideal of $I(A)$ \cite[Theorem 2]{halpern}. We aim to prove that $K=A\cap J$ is a liminal essential ideal
of $A$.

Suppose that $\alpha_0$ is an irreducible representation
of $J$ on a Hilbert space $H_{\alpha_0}$. As $J$ is an ideal of $I(A)$,
$\alpha_0$ extends uniquely to an irreducible representation $\alpha$ of $I(A)$
on the same Hilbert space $H_{\alpha_0}$.
Thus, $\alpha\left(I(A)\right)\supset\alpha(J)=\alpha_0(J)=K(H_{\alpha_0})$.

If $\hat J$ denotes the spectrum of $J$
(unitary equivalence classes of irreducible representations of $J$) and if,
for each $\alpha_0\in\hat J$, $\alpha$ denotes the unique
extension of $\alpha_0$ to an irreducible representation of
$I(A)$, we consider the representation $\rho$ of $I(A)$ defined by
\[
\rho\;=\;\bigoplus_{\alpha_o\in\hat J}\alpha\,.
\]
By construction, $\rho_{|J}$ is a faithful representation of $J$. We next show that $\rho_{|A}$ is a
faithful representation of $A$. Suppose that $a\in A^+$ satisfies $\rho(a)=0$.
If $e\in I(A)$ is any abelian projection, then $eae\in J$ and $\rho(eae)=\rho(e)\rho(a)\rho(e)=0$. Because $\rho_{|J}$
is a faithful representation of $J$, $eae=0$; so, $a^{1/2}e=0$.
Thus, $a^{1/2}e=0$ for all abelian projections of $I(A)$. Because $I(A)$ is a type I AW$^*$-algebra,
\[
1=\sup\,\{e\,:\,e\in I(A)\mbox{ is an abelian projection}\}\,.
\]
Therefore, by \cite[Lemma 1.9]{hamana1981},
\[
a\;=\;
a^{1/2}1a^{1/2}\;=\;\sup\,\{a^{1/2}ea^{1/2}\,:\,e\in I(A)\mbox{ is an abelian projection}\}\;=\;0\,,
\]
which proves that $\rho_{|A}$ is a faithful representation of $A$.

(Indeed $\rho$ is a faithful representation of $I(A)$ as well. To prove this, suppose that $h\in I(A)^+=\overline A^+$
satisfies $\rho(h)=0$. Thus, $\rho(a)=0$ for all $a\in A^+$ for which $a\le h$.
Since $\rho_{|A}$ is a faithful representation of $A$,
$\rho(a)=0$ only if $a=0$. Because $h=\sup\{a\in A^+\,:\,a\le h\}$ and $a=0$ for every $a\le h$, we
conclude that $h=0$, which proves that $\rho$ is faithful.)

Let $s\in J^+$ be nonzero and choose any $\alpha_0\in\hat J$.
Then $\alpha(a)$ is compact for every $a\in A^+$ such that $a\le s$. To verify this,
fix $a\in A^+$ for which $a\le s$; thus, $\alpha(a)\le\alpha(s)=\alpha_0(s)$.
Let $\{\xi_n\}_{n\in\mathbb N}$ be a sequence in the unit sphere of
the Hilbert space $H_{\alpha_0}$. By the compactness of $\alpha(s)^{1/2}$,
there is a subsequence $\{\xi_{n_k}\}_{k\in\mathbb N}$ such that
$\{\alpha(s)^{1/2}\xi_{n_k}\}_{k\in\mathbb N}$ is convergent. This implies that the sequence
$\{\alpha(a)^{1/2}\xi_{n_k}\}_{k\in\mathbb N}$
is a Cauchy sequence, for
\[
\begin{array}{rcl}
\|\alpha(a)^{1/2}\xi_{n_j}-\alpha(a)^{1/2}\xi_{n_m}\|^2\;&=&\;
\left\langle \alpha(a)\left(\xi_{n_j}-\xi_{n_m}\right),\,\left(\xi_{n_j}-\xi_{n_m}\right)\right\rangle \\
&&\\
\;&\le\;&
\left\langle \alpha(s)\left(\xi_{n_j}-\xi_{n_m}\right),\,\left(\xi_{n_j}-\xi_{n_m}\right)\right\rangle \\
&&\\
\;&=&\;
\|\alpha(s)^{1/2}\xi_{n_j}-\alpha(a)^{1/2}\xi_{n_m}\|^2\,.
\end{array}
\]
Hence, $\{\alpha(a)^{1/2}\xi_{n_k}\}_{k\in\mathbb N}$ is convergent, which yields $\alpha(a)$ compact.
Since the choice of $\alpha_0\in\hat J$ is arbitrary, this shows that $\rho(a)\in \rho(J)$ if $a\in A^+$
satisfies $a\le s$. Because $\rho$ is faithful, this is to say that $a\in J$ if
$a\in A^+$ satisfies $a\le s$. Furthermore, since $s$ is nonzero and $A$ is order dense in $I(A)$,
there is a nonzero $a\in A^+$ such that $a\le s$. In particular, this nonzero $a$ belongs to $J$,
thereby proving that $K=A\cap J\not=\{0\}$.

The previous paragraph establishes the following identity:
\[
s\;=\;\sup\,\{a\in K^+\,:\,a\le s\}\,,\quad\forall\,s\in J^+\,.
\]
This fact will now be used to prove that $K$ is an essential ideal of $A$.
To this end, let $L$ be any ideal of $A$ for which $L\cap K=\{0\}$. Thus if $b\in L^+$, then $bab=0$
for every $a\in K^+$. Now, if $e\in I(A)$ is any abelian projection, then $e\in J^+$ and
\[
e\;=\;\sup\,\{a\in K^+\,:\,a\le e\}\,.
\]
Therefore, again by \cite[Lemma 1.9]{hamana1981},
\[
beb\;=\;\sup\,\{bab\in K^+\,:\,a\le e\}\;=\;0\,.
\]
Thus, $eb=be=0$ for every abelian projection $e\in I(A)$, which implies that $b=0$ (as demonstrated
earlier in this proof). Hence, $L\cap K=\{0\}$ only if $L=\{0\}$ and so $K$ is an essential ideal of $A$.

The final point to verify is that $K$ is liminal. But this follows from the fact that every C$^*$-subalgebra
of a liminal C$^*$-algebra is liminal \cite[Proposition 4.2.4]{dixmier}, and by noting that $K$ is a C$^*$-subalgebra
of the liminal ideal $J$ of $I(A)$.
\qed\medskip

\section{Applications}

\begin{theorem}\label{applics} The following statements hold for every separable C$^*$-algebra $A$.
\begin{enumerate}
\item\label{lim-postlim} $A$ has a liminal essential ideal if and only if $A$ has postliminal essential ideal.
\item \label{tomato} If $A$ is abelian, then $I(A)$ is a W$^*$-algebra if and only if
there exists a finite or countably infinite set $\Gamma$ such that $I(A)=l^\infty(\Gamma)$ and
$c_0(\Gamma)\subseteq A\subseteq l^\infty(\Gamma)$.
\item\label{onion} If $A$ is simple and $I(A)$ is a W$^*$-algebra, then $A=K(H)$ for some Hilbert space $H$.
\item\label{faithful state} $I(A)$ admits a faithful state.
\item\label{prime2} If $A$ is prime, then exactly one of the following two
statements holds:
\begin{enumerate}
\item $I(A)\cong_* B(H)$, for some separable Hilbert space $H$;
\item $I(A)$ is a wild type III AW$^*$-factor.
\end{enumerate}
In particular, if $A$ has no postliminal
essential ideal, then $I(A)$ is a wild type III AW$^*$-factor.
\end{enumerate}
\end{theorem}

\proof For the proof of (\ref{lim-postlim}), every liminal ideal is postliminal, by definition. Thus,
assume that $A$ has a postliminal essential ideal, say $K$.
As $A$ and $K$ are separable and $K$ is an essential ideal,
$\overline K=\overline A$ (Theorem \ref{sigma completions}).
Because $K$ is liminal, $\overline K$ is type I, and so $\overline A=I(A)$ is of type I.
By Theorem \ref{aw-typeI}, $A$  has a liminal essential ideal, which proves (\ref{lim-postlim}).

To prove (\ref{tomato}), suppose now that $A$ is abelian and $I(A)$ is
a W$^*$-algebra. By Theorem \ref{envelope is type I},
$A$ has a minimal essential ideal $K$ for which $K\cong_*\bigoplus_{n\in\Gamma} K(H_n)$,
for some finite or countable infinite set $\Gamma$; however, as $K$ is
abelian, $K(H_n)=\mathbb C$ for every $n$, whence $K=c_0(\Gamma)$. As $K$ contains all minimal
projections in $A$, we deduce that $A\subseteq l^\infty(\Gamma)$. Finally, $I(c_0(\Gamma))=l^\infty(\Gamma)$
(because $c_0(\Gamma)$ is order dense in $l^\infty(\Gamma)$),
so that $I(A)=l^\infty(\Gamma)$.
The converse is a direct application of Theorem \ref{envelope is type I} where the minimal essential ideal
of $A$ is $c_0(\Gamma)$.

To prove (\ref{onion}),
assume that $A$ is simple and that $I(A)$ is a W$^*$-algebra.
By (\ref{quatre}) of Theorem \ref{envelope is type I}, $A$ has a minimal
essential ideal of the form $K=\bigoplus_nK(H_n)$. Being simple, $A=K$; and for $K$ to be simple,
there can be only one summand. Thus, $A=K(H_1)$.

For the proof of (\ref{faithful state}), note that because $A$ is separable, $A$ has a faithful representation
as a C$^*$-subalgebra of $B(H)$,
where $H$ is a separable Hilbert space. Thus, by Hamana's construction of the injective envelope,
there is a projection $\phi:B(H)\rightarrow B(H)$ such that $\phi\left(B(H)\right)=I(A)$. The separability of $H$
implies that $B(H)$ has a faithful state $\omega$. This state is also faithful on the C$^*$-algebra representation
of $I(A)$. To prove this, recall that the product $\circ$ on $I(A)$ is given by $x\circ y=\phi(xy)$, for all
$x,y\in I(A)$. Suppose $x\in I(A)$ is such that $\omega(x^*\circ x)=0$. Then $\omega\left(\phi(x^*x)\right)=0$ and so
$\phi(x^*x)=0$, as $\omega$ is a faithful state on $B(H)$. Therefore, by the Schwarz inequality for completely
positive maps, $0\le\phi(x)^*\phi(x)\le\phi(x^*x)=0$. This implies that $\phi(x)=0$. However, on $I(A)$ the map
$\phi$ acts as the identity. Thus, $x=\phi(x)=0$, which proves that $\omega$ is a faithful state on the C$^*$-algebra
representation of $I(A)$.

To prove (\ref{prime2}), assume now that $A$ is prime. By Proposition \ref{factors}, $I(A)$ is a factor.
But this factor cannot be of type II for the following reasons.
Proposition \ref{no type II_1} already excludes the case of finite type II AW$^*$-factors.
By \cite{elliott--saito--wright},
every type II${}_\infty$ AW$^*$-factor that admits a faithful state is a W$^*$-factor.
Since $I(A)$ admits a faithful state and
since $I(A)$ is a W$^*$-algebra only in the case where $I(A)$
is of type I (Theorem \ref{envelope is type I}), it is
impossible for $I(A)$ to be a type II${}_\infty$ AW$^*$-factor.
Hence, $I(A)$ is a factor of either type I or type III.

In the case where $I(A)$ is of type I we have $I(A)\cong_*B(H)$ for some Hilbert space $H$,
because all type I AW$^*$-factors have this form \cite[Theorem 2]{kaplansky1952}. Indeed, in this case,
$\overline A\,{}^\sigma=I(A)\cong_* B(H)$; since $\overline A\,{}^\sigma$ is countably decomposable,
$H$ can be chosen to be separable.

If $I(A)$ is not of type I, then the type III AW$^*$-factor $I(A)$ is cannot be a W$^*$-algebra, by
Theorem \ref{envelope is type I}. Every AW$^*$-factor that is not W$^*$-algebra is wild \cite{wright1976b}; hence,
$I(A)$ is wild.
\qed\medskip

We wish to remark that statement (\ref{faithful state}) of Theorem \ref{applics} above was previously
noted (without proof) and employed
in \cite[Corollary 3.8]{hamana1981}.
\medskip

Turning now to the local multiplier algebra, in most
cases the precise determination of $M_{\rm loc}(A)$ is difficult,
and so one is interested to know what properties $M_{\rm loc}(A)$ might exhibit.
In particular, the following
questions have been raised in the literature.
\begin{enumerate}
\item[(Q1)] For which C$^*$-algebras $A$ is $M_{\rm loc}\left(M_{\rm loc}(A)\right)= M_{\rm loc}(A)$\,?
(\cite{ara--mathieu2003,somerset2000})
\item[(Q2)] For which C$^*$-algebras $A$ is $M_{\rm loc}(A)$ injective\,? (\cite{frank2002,frank--paulsen})
\end{enumerate}
Partial answers to these questions are listed in the
theorem below.

\begin{theorem} Assume that $A$ is a separable C$^*$-algebra.
\begin{enumerate}
\item\label{lettuce} If $A$ has a liminal essential ideal,
then $M_{\rm loc}\left(M_{\rm loc}(A)\right)$ is an injective
C$^*$-algebra of type I and
\[
M_{\rm loc}\left(M_{\rm loc}(A)\right)\;=\;
\overline A\;=\;I(A)\,.
\]
\item\label{egg} If $A$ has a minimal essential ideal that is $*$-isomorphic to a C$^*$-algebraic direct sum
of elementary C$^*$-algebras, then $M_{\rm loc}(A)$ is an injective W$^*$-algebra of type I and
\[
M_{\rm loc}(A)\;=\;M_{\rm loc}\left(M_{\rm loc}(A)\right)\;=\;\overline A\;=\;I(A)\,;
\]
\end{enumerate}
\end{theorem}

\proof To prove (\ref{lettuce}), let $K$ be a liminal essential ideal of $A$.
As $A$ and $K$ are separable and $K$ is an essential ideal,  $\overline K=\overline A$.
Because $K$ is liminal, $\overline K$ is type I, and so $\overline A=I(A)$ is of type I.
Again using that $A$ and $K$ are separable and that $K$ is an essential ideal, conclude that from
\cite[Theorem 2.8]{somerset2000} that
$M_{\rm loc}\left(M_{\rm loc}(A)\right)=\overline A$.
Hence,
$M_{\rm loc}\left(M_{\rm loc}(A)\right)$ is an injective
C$^*$-algebra of type I.

For the proof of (\ref{egg}), note that
Theorem \ref{envelope is type I} and its proof imply there is a minimal essential ideal $K$ of $A$ such that
$K\cong_*\bigoplus K(H_n)$ and $M(K)=M_{\rm loc}(A)=\overline A=I(A)$.
Every
AW$^*$-algebra is its own local multiplier algebra \cite[Theorem 2.3.8]{ara--mathieu2003}, and so
\[
M_{\rm loc}\,(A)\,=\,M_{\rm loc}\left(M_{\rm loc}(A)\right)\,=\, I(A)\,.
\]
This completes the proof of (\ref{egg}).
\qed\medskip

There is an unresolved issue: is $M_{\rm loc}(A)$ injective if $A$ is separable and has a liminal
essential ideal\,? Recall that if $K$ is an essential ideal of $A$,
then $\overline K=\overline A$
\cite{saito1980}. Thus, it is sufficient to ask: is $M_{\rm loc}(A)$ injective if $A$ is separable and liminal\,?
This question is at present open.

\section{Nonseparable C$^*$-algebras}

The focus of this paper has been on separable C$^*$-algebras. For example,
Proposition \ref{no type II_1} and Theorem \ref{envelope is type I}
do not hold for nonseparable C$^*$-algebras. More
specifically, if $R$ denotes the
hyperfinite II${}_1$ factor $R$, then $R$ is injective and, thus, $R=I(R)$ is
a W$^*$-factor of type II. However, this leads to another question of interest:
if $M$ is a nonhyperfinite II${}_1$ factor,
then what is the injective envelope of $M$\,? Because $M$ is simple, $I(M)$
is an AW$^*$-algebra factor; is $I(M)$ a
finite AW$^*$-factor\,? More generally, does the passage from $M$ to $I(M)$
preserve type if $M$ is a von Neumann algebra\,?

Although Proposition \ref{no type II_1} and Theorem \ref{envelope is type I} do not hold for nonseparable
C$^*$-algebras, the necessity part of Theorem \ref{aw-typeI} was established without recourse to separability.
Thus, the following theorem holds.

\begin{theorem} If the injective envelope of a
C$^*$-algebra $A$ is of type I, then $A$ has a liminal
essential ideal.
\end{theorem}


The original motivation for the concept of injectivity is
Arveson's Hahn--Banach Extension Theorem \cite{arveson1969}
for completely positive linear maps, and the idea of an injective
envelope stems from Arveson's theory of
boundary representations \cite{arveson1972}. In the work on boundary representations,
the algebras under consideration need
not have been separable, but frequently the algebras were assumed to have nontrivial
intersection with the compact operators.
In this spirit we have the following result, which generalises one form
the ``boundary theorem" from $B(H)$ to discrete type I von Neumann algebras and which shows that
statement (\ref{trois}) of Theorem \ref{envelope is type I} holds for nonseparable C$^*$-algebras as well.

\begin{theorem}\label{atomic} If
$\pi:A\rightarrow B(H)$ is a faithful representation of a C$^*$-algebra $A$ on a Hilbert space $H$
such that $\pi(A)''$ is generated by its minimal projections, each of which is contained in
$\pi(A)$, then $\pi(A)''=I(A)$.
\end{theorem}
\proof
Without loss of generality,
we may assume that $A$ is already faithfully represented as a C$^*$-subalgebra
of $B(H)$ such that $M=A''$ is generated by its minimal projections, each of which is contained in $A$.
Because $M$ is generated by minimal projections, $M$ is an injective von Neumann algebra,
by Corollary \ref{atomic-injective}.
To show that $M$ is the injective envelope of $A$, it is sufficient, by Proposition \ref{hamana37}, to show that any
completely positive linear map $\varphi:M\rightarrow M$ that fixes $A$ must be the identity map on $M$.
If this is indeed so, then $M$ is an injective envelope for $A$ and, by the uniqueness of
the injective envelope, we deduce that $M=I(A)$.
If $\varphi:M\rightarrow M$ is a completely positive map
such that $\varphi_{|A}=\mbox {id}_A$, then we will show that $\varphi=\mbox{id}_M$.

To this end, observe that because $\varphi:M\rightarrow M$ is a unital completely positive map that preserves
$A$, $\varphi$ has the following property:
\[\varphi(xk)=\varphi(x)k,\mbox{ for every }k\in A.\]
This fact follows from the Cauchy-Schwarz
inequality and from the fact that $A$ is in the multiplicative domain of $\varphi$ (see \cite[9.2]{stratila} or
\cite[Corollary 2.6]{ozawa2004}).
Using this fact
we shall
deduce below that
\begin{equation}\label{reversa}x\ge0\mbox{ if and only if
}\varphi(x)\ge0.\end{equation}
Indeed, one implication is obvious from the positivity of $\varphi$.
To prove the other implication,
assume that $\varphi(x)\ge0$. Thus, $\varphi(\im(x))=\im(\varphi(x))=0$.
Let $z=\im(x)$ and write $z=z^+-z^-$, where $z^+,z^-\in M^+$ are such that $z^+z^-=z^-z^+=0$.

Our first goal is to prove that $z^+=0$. Suppose, on the contrary, that $z^+\not=0$. Thus, there is
a strictly positive $\lambda$ in the spectrum of $z^+$; hence,
there is a spectral projection $p\in M$ such that $0\not=\lambda p\le pz^+=z^+p$.
Note that $z^-p=0$, as the projection $p$ is in the von Neumann algebra generated by $z^+$ and $z^+z^-=z^-z^+=0$.
Let $q\in A$ be an arbitrary
minimal projection of $M$ and consider the projection $p\wedge q\in M$. Because $p\wedge q\le q$
and $q$ is minimal, either $p\wedge q=0$ or $p\wedge q=q$. We will show that
the latter case cannot occur (under the conventional assumption that minimal projections are defined
to be nonzero). Assume that it is true that $p\wedge q=q$. Then $0\not=q=p\wedge q\le p$.
Pre- and post-multiply
the inequality $\lambda q\le \lambda p\le z^+p=zp$ by $q$ to obtain $\lambda q\le q(zp)q\le qzq$.
Note that $\varphi(zq)=\varphi(z)q$ (because $A$
is in the multiplicative domain of $\varphi$) and that $\varphi(z)=0$ (by hypothesis).
Likewise, for any hermitian $y\in M$,
$\varphi(qy)=\varphi(yq)^*=q\varphi(y)$. Thus, $\varphi(qzq)=q\varphi(z)q=0$ and
$0\le\lambda q=\varphi(\lambda q)\le q\varphi(z)q=0$. This implies that $q=0$, which contradicts the fact
that $q$ is minimal and, thus, nonzero. Therefore, it must be that $p\wedge q=0$, for every minimal projection
$q$ of $M$. Because every nonzero projection in $M$
majorises a minimal projection, we conclude that $p=0$, in contradiction
to the fact that $p$ is a nonzero spectral projection of $z^+$. Hence, it must be that
$z^+=0$.

A similar argument shows that $z^-=0$.
We can find a nonzero $\lambda\in\mathbb R^+$ and a minimal projection $q\in A$ such that $qzq\le-\lambda q$; thus
$-\lambda q=\varphi(-\lambda q)\ge\varphi(qzq)=q\varphi(z)q=0$, and again $q=0$.

We conclude that $z=0$, which implies that
$x$ is selfadjoint. It remains to show that $x$ is positive.
Assume that $x$ is not positive.
Thus, there exists a nonzero spectral projection in the
negative part of $\sigma(x)$; by taking once again a suitable minimal
subprojection $q$, we can find $\lambda>0$ such that
$qxq\le-\lambda q$. But then $\varphi(qxq)\le-\lambda q$; and on the
other hand, $\varphi(qxq)=q\varphi(x)q\ge0$. The
contradiction implies that no such $q$ can exist, and so $x\ge0$.

From (\ref{reversa}) and the fact the $\varphi$ preserves $A$,
we have that $k\in A$, $k\le x$ if
and only if $k\le\varphi(x)$. Statement (\ref{A has min projs}) of Lemma \ref{prop-sup}
asserts that $A$ is order dense in $M$. Hence,
$\varphi(x)=x$ for every $x\in M^+$, which implies that
$\varphi$ is the identity map on $M$.
\qed\medskip

\section{Open Questions}

Although this paper is mainly concerned with type I injective envelopes of separable
C$^*$-algebras, there are a number of unresolved questions that underscore the limits of
our current state of knowledge concerning injective envelopes in general. A few such questions are listed here.

\begin{enumerate}
\item Suppose that $A$ is a separable C$^*$-algebra.
\begin{enumerate}
\item Is $M_{\rm loc}(A)$ an AW$^*$-algebra\,?
\item Is $\overline A=I(A)$ if $I(A)$ is of type III\,?
\end{enumerate}
\item Suppose that $\displaystyle\bigotimes_1^\infty M_2$ and
$\displaystyle\bigotimes_1^\infty M_3$ denote the UHF C$^*$-algebras obtained
through the tensor products of the matrix algebras $M_2$ and $M_3$ respectively.
The injective envelope of each of these C$^*$-algebras is a wild type III AW$^*$-factor. Is it
true that
\[
I\left(\displaystyle\bigotimes_1^\infty M_2\right)\;=\;I\left(\displaystyle\bigotimes_1^\infty M_3\right)\,?
\]
\item Suppose that $M$ is a von Neumann algebra.
\begin{enumerate}
\item What is $I(M)$ if $M$ is not injective\,?
\item If $M$ is a non-injective type II${}_1$ factor, then is the AW$^*$-factor $I(M)$ also of type II\,?
\end{enumerate}
\end{enumerate}



\begin{thebibliography}{99}

\bibitem{ara--mathieu2003}
            P.~Ara and M.~Mathieu,
            {\em Local Multipliers of C$^*$-algebras},
            Springer Monographs in Mathematics,
            London, 2003.

\bibitem{arveson1969}
           W.B.~Arveson,
            Subalgebras of C$^*$-algebras,
            {\em Acta Math.}
            123 (1969), 141--224.

\bibitem{arveson1972}
            W.B.~Arveson,
            Subalgebras of C$^*$-algebras, II,
            {\em Acta Math.}
            128 (1972), 271--308.

\bibitem{choi1974}
            M.-D.~Choi,
            A Schwarz inequality for positive linear maps on C$^*$-algebras,
            {\em Illinois J.~Math.}
            18 (1974), 565--574.

\bibitem{choi--effros1977}
            M.-D.~Choi and E.G.~Effros,
            Injectivity and operator spaces,
            {\em J.~Functional Analysis}
            24 (1977), 156--209.

\bibitem{dixmier}
            J.~Dixmier,
            {\em Les C$\,{}^*$-alg\`ebres et leurs repr\'esentations},
            Gauthier--Villars,
            Paris, 1969.

\bibitem{elliott--saito--wright}
            G.A.~Elliott, K.~Sait\^o, and J.D.M.~Wright,
            Embedding AW$^*$-algebras as double commutants in type I algebras,
            {\em J.~London Math.~Soc.}
            28 (1983), 376--384.

\bibitem{effros--ruan}
            E.G.~Effros and Z.-J.~Ruan,
            {\em Operator Spaces},
            Oxford University Press, Oxford, 2000.

\bibitem{frank2002}
            M.~Frank,
            Injective envelopes and local multiplier algebras of C$^*$-algebras,
            {\em Int.~Math.~J.}
            1 (2002), 611--620.

\bibitem{frank--paulsen}
            M.~Frank and V.I.~Paulsen,
            Injective envelopes of C$^*$-algebras as operator modules,
            {\em Pacific J.~Math.}
            212 (2003), 57--69.

\bibitem{hamana1979}
            M.~Hamana,
            Injective envelopes of C$^*$-algebras,
            {\em J.~Math.~Soc.~Japan}
            31 (1979), 181--197.

\bibitem{hamana1981}
            M.~Hamana,
            Regular embeddings of C$^*$-algebras in monotone complete C$^*$-algebras,
            {\em J.~Math.~Soc.~Japan}
            33 (1981), 159--183.

\bibitem{halpern}
            H. Halpern,
            The maximal GCR-ideal in an AW$^*$-algebra,
            {\em Proc. Amer. Math. Soc.}
            17 (1966), 906--914.

\bibitem{kadison1956}
            R.V.~Kadison,
            Operator algebras with a faithful weakly-closed representation,
            {\em Ann. of Math.}
            64 (1956), 175--181.

\bibitem{kaplansky1952}
            I.~Kaplansky,
            Algebras of type I,
            {\em Ann. of Math.}
            56 (1952), 460--472.

\bibitem{ozawa2004}
            N.~Ozawa,
            About the QWEP conjecture,
            {\em International J.~Math.}
            15 (2004), 501--530.

\bibitem{ozawa--saito}
            M.~Ozawa and K.~Sait\^o,
            Embeddable AW$^*$-algebras and regular completions,
            {\em J.~London Math.~Soc.}
            34 (1986), 511--523.

\bibitem{paulsen}
            V.I.~Paulsen,
            {\em Completely Bounded Maps and Operator Algebras},
            Cambridge University Press, Cambridge, 2002.

\bibitem{pedersen1979}
            G.K.~Pedersen,
            {\em C$^*$-Algebras and their Automorphism Groups},
            Academic Press, London, 1979.

\bibitem{saito1980}
            K.~Sait\^o,
            A structure theory in the regular $\sigma$-completion of a C$^*$-algebra,
            {\em J.~London Math.~Soc.}
            22 (1980), 549--558.

\bibitem{somerset2000}
            D.W.B.~Somerset,
            The local multiplier algebra of a C$^*$-algebra, II,
            {\em J.~Functional Analysis}
            171 (2000), 308--330.

\bibitem{stratila}
            \c{S}. Str\u{a}til\u{a},
            {\em Modular Theory in Operator Algebras},
            Editura Academiei, Bucharest, 1981.

\bibitem{tak}
            M.~Takesaki,
            {\em Theory of Operator Algebras I},
            Encylopaedia of Mathematical Sciences, Vol 124,
            Springer-Verlag, New York, 2001.

\bibitem{tomiyama}
            J.~Tomiyama,
            Tensor products and projections of norm one in von Neumann algebras,
            {\em Mimeographed Lecture Notes},
            Mathematical Institute, Copenhagen University, Copenhagen, 1970.

\bibitem{wright1976a}
            J.D.M.~Wright,
            Regular $\sigma$-completions of C$^*$-algebras,
            {\em J.~London Math.~Soc.}
            12 (1976), 299--309.


\bibitem{wright1976b}
            J.D.M.~Wright,
            Wild AW$^*$-factors and Kaplansky--Rickart algebras,
            {\em J.~London Math.~Soc.}
            13 (1976), 83--89.

\end{thebibliography}
\end{document}